\documentclass[a4paper,12pt]{article}
\usepackage[paper=letterpaper,margin=1in]{geometry}
\usepackage{graphicx}
\usepackage{hyperref}
\hypersetup{
    colorlinks,
    citecolor=red,
    filecolor=black,
    linkcolor=blue,
    urlcolor=black
}

\usepackage{amsmath,amssymb,amsfonts,epsfig,cite,setspace,bigstrut, verbatim, multirow}
\pdfoutput=1

\makeatother
\date{}
\numberwithin{equation}{section}

\begin{document}

\title{\textbf{Hecke Groups, \emph{Dessins d'Enfants} and the Archimedean Solids}}

\author{Yang-Hui He\textsuperscript{1,2}\footnote{hey@maths.ox.ac.uk}~ and
James Read\textsuperscript{1}\footnote{james.read@merton.ox.ac.uk}~
}
\begin{onehalfspace}
\maketitle
\begin{center}
\textsuperscript{1}\emph{Merton College, University of Oxford, OX1 4JD, UK}\\
~\\
\textsuperscript{2}\emph{Department of Mathematics, City University,
London,}\\
\emph{Northampton Square, London EC1V 0HB, UK;}\\
\emph{School of Physics, NanKai University, Tianjin, 300071, P.R.
China;}\\
\par\end{center}

\medskip{}

\begin{abstract}

Grothendieck's \emph{dessins d'enfants} arise with ever-increasing frequency in many areas of 21\textsuperscript{st} century mathematical physics.
  In this paper, we review the connections between \emph{dessins} and the theory of \emph{Hecke groups}. Focussing on the restricted class of highly symmetric \emph{dessins} corresponding to the so-called \emph{Archimedean solids}, we apply this theory in order to provide a means of computing representatives of the associated conjugacy classes of Hecke subgroups in each case.
The aim of this paper is to demonstrate that \emph{dessins} arising in mathematical physics can point to new and hitherto unexpected directions for further research. In addition, given the particular ubiquity of many of the \emph{dessins} corresponding to the Archimedean solids, the hope is that the computational results of this paper will prove useful in the further study of these objects in mathematical physics contexts.

\end{abstract}
\end{onehalfspace}

\pagebreak{}

\tableofcontents{}


\section{Introduction}

Grothendieck's \emph{dessins d'enfants} -- bipartite graphs drawn on Riemann surfaces -- arise with ever-increasing frequency in 21\textsuperscript{st} century mathematical physics, appearing in e.g.~the study of $\mathcal{N}=2$ supersymmetric gauge theories \cite{FirstCachazo,CachazoDessins,HR2015}, elliptically fibred Calabi-Yaus \cite{YMR,1,He:2013lha}, toric conformal field theories \cite{T1,T2}, topological strings \cite{S1, S2}, and elsewhere. Given this growing ubiquity, it is valuable to study both the basic theory of \emph{dessins}, and the applications of this theory to particularly significant cases. This paper aims to achieve both such goals, by first reviewing the mathematical results which connect \emph{dessins d'enfants} to the theory of \emph{Hecke groups}, before applying this work to the \emph{dessins} for the so-called \emph{Archimedean solids}. The hope is to provide a compendium of computational results to aid future research featuring these objects.

With the above in mind, let us begin our story by recalling that
the {\em Platonic solids} are  the convex polyhedra with equivalent faces composed of
 congruent convex regular polygons; these objects have been known and studied
  for millennia. These solids also appear in mathematical physics:~for example, the symmetries of the Platonic solids are related to D-brane orbifold theories \cite{HN}. Now, more broadly, another well-known class of convex
 polyhedra are the {\em Archimedean solids}:~the semi-regular convex polyhedra composed of two or more types of regular polygons meeting in
 identical vertices, with no requirement that faces be equivalent.
There are three categories of such Archimedean solids:~(I) the Platonic solids; (II) two infinite series solutions -- the prisms and
anti-prisms; and (III) fourteen further exceptional cases.

In \cite{Zvon}, a novel approach to the Archimedean solids was taken by interpreting the graphs
of these solids as clean {\em dessins d'enfants}. By \emph{clean}, we mean that
all the nodes of one of the two possible colours of the bipartite graph have valency two \cite{4,DessinBook}. 
Now, the planar graph of any polytope can be interpreted as a clean \emph{dessin} by inserting a black node into every edge of the graph, and colouring every vertex white. In this way, one can construct \emph{dessins} for all the Archimedean solids. In each case, the underlying Riemann surface is the sphere $\mathbb{C}\mathbb{P}^1$, and the \emph{dessins} are drawn in a planar projection.

 
In a parallel vein, \emph{trivalent} clean \emph{dessins} can be associated to conjugacy classes of subgroups
 of the modular group $\Gamma = \mathrm{PSL}\left(2,\mathbb{Z}\right)$. (Recall that the modular group $\Gamma\equiv\Gamma\left(1\right)$ is the group of linear fractional transformations $z\rightarrow\frac{az+b}{cz+d}$, 
with $a,b,c,d\in\mathbb{Z}$ and $ad-bc=1$. The presentation of
$\Gamma$ is $\left\langle x,y|x^{2}=y^{3}=I\right\rangle$.) To do so, we replace all $n$-valent vertices of the \emph{dessin} with oriented $n$-gons. This constructs a \emph{Schreier coset graph}, which displays the permutation action of generators $x$ and $y$ on each coset $G x_i$ of $\mathrm{PSL}\left(2,\mathbb{Z}\right)=\bigcup_i G x_i$, $i=1, \ldots, \mu$, where $\mu$ is the index of some subgroup $G$ in $\mathrm{PSL}\left(2,\mathbb{Z}\right)$ \cite{2}.
 Generalising, \emph{any} clean \emph{dessin} (not necessarily trivalent) 
  can be associated with a conjugacy class of subgroups of a certain so-called \emph{Hecke group} $H_{n}$, defined as having presentation $\left\langle x,y|x^{2}=y^{n}=I\right\rangle$.


Once a Hecke subgroup $G$ has been associated to the \emph{dessin} in question, further results follow. For example, one can quotient the upper half plane by $G$ to construct the surface $\mathcal{H}/G$. From this, one can construct a Belyi map (in a manner detailed in the main body of this paper), i.e.~a holomorphic map to $\mathbb{P}^1$ \emph{ramified} at only $\left\{ 0, 1, \infty \right\}$. To each Belyi map there corresponds a unique \emph{dessin}; precisely the \emph{dessin} with which we began! Hence, this chain of connections compactifies to a circle.

Returning now to the Archimedean solids, we see that, interpreting these as \emph{dessins}, we should be able to explore the above circle of connections via explicit computations. Indeed, the Belyi maps associated to these \emph{dessins} have already been computed in \cite{Zvon}; therefore, our task is to complete the circle by computing representatives of the conjugacy classes of Hecke subgroups associated to these objects. Not only will this provide an illustration of this aspect of the mathematical theory underlying \emph{dessins}, but it will also provide a useful resource for further research in this area. Indeed, given that some \emph{dessins} for the Archimedean solids have have already arisen in areas of mathematical physics (see for example \cite{YMR,1,HR2015}), it is not unreasonable to expect that such \emph{dessins} will continue to manifest themselves in future research.

The structure of this paper is as follows.
In section 2 we present some technical details regarding Hecke groups and \emph{dessins d'enfants}.
We show that it is possible to interpret each clean \emph{dessin} as the {Schreier coset graph}
for a conjugacy class of subgroups of a certain Hecke group, and review the circle of connections mentioned above.
 In section 3, we
find the permutation representations for the conjugacy classes of subgroups of Hecke groups corresponding
to every Archimedean solid, and provide an algorithm to compute explicit generating sets of matrices
for representatives of these conjugacy classes
in each case. In section 4, we close with some conclusions, returning to the specific applications of this work in various subfields of mathematical physics.

\section{Dessins d'Enfants and Hecke Groups}

In this section, we review some essential details regarding both Hecke groups and clean \emph{dessins
d'enfants}. We begin by considering the modular group $\Gamma\cong H_{3}$, as
although this is isomorphic to only one particular Hecke group, it is by far the most well-studied, and we shall
draw upon the presented results at several points in the ensuing discussion. Subsequently, we discuss
Hecke groups more generally, before moving on to consider clean  \emph{dessins} and their associated Belyi maps.
Note that the connection between
trivalent  \emph{dessins} and the modular group is discussed in \cite{JS1, JS2}, while the relation between Hecke groups and maps
is discussed in e.g. \cite{normal, Jones}.
With these results in hand, we describe how every clean  \emph{dessin} is isomorphic to the
Schreier coset graph for a conjugacy class of subgroups of a Hecke group, and how the circle of connections closes via a correspondence between such subgroups and the Belyi maps for the original  \emph{dessins}.

\subsection{The Modular Group}

To begin, let us recall some essential details regarding the modular group $\Gamma$.
This is the group of linear fractional transformations $\mathbb{Z} \ni z\rightarrow\frac{az+b}{cz+d}$, 
with $a,b,c,d\in\mathbb{Z}$ and $ad-bc=1$. It is generated by the
transformations $T$ and $S$ defined by:
\begin{equation}
T(z)=z+1\quad,\quad S(z)=-1/z \ .
\end{equation}
The presentation of $\Gamma$ is $\left\langle S,T|S^{2}=\left(ST\right)^{3}=I\right\rangle $, and we will later discuss the presentations of certain modular subgroups. The
$2\times2$
matrices for $S$ and $T$ are as follows:
\begin{equation}
T=\begin{pmatrix}1 & 1\\
0 & 1
\end{pmatrix}\quad,\quad S=\begin{pmatrix}0 & -1\\
1 & 0
\end{pmatrix} \ .
\end{equation}

Letting $x = S$ and $y = ST$ denote the elements of order 2 and 3 respectively,
we see that $\Gamma$ is the free product of the cyclic groups $C_{2}=\left\langle x|x^{2}=I\right\rangle $
and $C_{3}=\left\langle y|y^{3}=I\right\rangle $. 
It follows that $2\times2$ matrices for $x$ and $y$ are:
\begin{equation}\label{xy}
x=\begin{pmatrix}0 & -1\\
1 & 0
\end{pmatrix}\quad,\quad y=\begin{pmatrix}0 & -1\\
1 & 1
\end{pmatrix} \ .
\end{equation}

With these details in hand, we can consider some important subgroups of $\Gamma$.

\subsubsection{Congruence Modular Subgroups}

The most significant subgroups of $\Gamma$ are the \emph{congruence} subgroups, defined by having the the entries in the generating matrices $S$ and $T$ obeying some modular arithmetic. Some conjugacy classes of congruence subgroups of particular
note are the following:

\begin{itemize}
\item Principal congruence subgroups:
\[
\Gamma\left(m\right):=\left\{ A\in\mathrm{SL}(2;\mathbb{Z})\;;\; A\equiv\pm I\;\mathrm{mod}\; m\right\} /\left\{ \pm I\right\} ;
\]

\item Congruence subgroups of level $m$: subgroups of $\Gamma$ containing
$\Gamma\left(m\right)$ but not any $\Gamma\left(n\right)$ for $n<m$;
\item Unipotent matrices:
\[
\Gamma_{1}\left(m\right):=\left\{ A\in\mathrm{SL}(2;\mathbb{Z})\;;\; A\equiv\pm\begin{pmatrix}1 & b\\
0 & 1
\end{pmatrix}\;\mathrm{mod}\; m\right\} /\left\{ \pm I\right\} ;
\]

\item Upper triangular matrices:
\[
\Gamma_{0}\left(m\right):=\left\{ \begin{pmatrix}a & b\\
c & d
\end{pmatrix}\in\Gamma\;;\; c\equiv0\;\mathrm{mod}\; m\right\} /\left\{ \pm I\right\} .
\]

\item The congruence subgroups
\[
\Gamma\left(m;\frac{m}{d},\epsilon,\chi\right):=\left\{ \pm\begin{pmatrix}1+\frac{m}{\epsilon\chi}\alpha & d\beta\\
\frac{m}{\chi}\gamma & 1+\frac{m}{\epsilon\chi}\delta
\end{pmatrix}\;;\;\gamma\equiv\alpha\;\mathrm{mod}\;\chi\right\} 
\]
for certain choices of $m$, $d$, $\epsilon$, $\chi$ (see \cite{1,YMR}).

\end{itemize}

We note here that:
\begin{equation}
\Gamma\left(m\right)\subseteq\Gamma_{1}\left(m\right)\subseteq\Gamma_{0}\left(m\right)\subseteq\Gamma \ .
\end{equation}

In section 3 of this paper, we shall remark on
the connections between some specific conjugacy classes of congruence modular subgroups and the
Archimedean solids.

\subsection{Hecke Groups}\label{s:hecke}

We can now extend our discussion of the modular group $\Gamma\cong H_{3}$ to the more
general Hecke groups $H_{n}$.
The Hecke group $H_{n}$ has presentation $\left\langle x,y|x^{2}=y^{n}=I\right\rangle$,
and is thus the free product of cyclic groups $C_{2}=\left\langle x|x^{2}=I\right\rangle$
and $C_{n}=\left\langle y|y^{n}=I\right\rangle$. 
Note that $\Gamma\cong H_{3}$
where $\Gamma$ is the modular group; and that $H_n$ is the triangle group $\left(2,n,\infty\right)$.
$H_{n}$ is generated by
transformations $T$ and $S$ now defined by:
\begin{equation}
T(z)=z+\lambda_{n}\quad,\quad S(z)=-1/z \ ,
\end{equation}
where $\lambda_{n}$ is some real number to be determined. The $2\times2$
matrices for these $S$ and $T$ are:
\begin{equation}
T=\begin{pmatrix}1 & \lambda_{n}\\
0 & 1
\end{pmatrix}\quad,\quad S=\begin{pmatrix}0 & -1\\
1 & 0
\end{pmatrix} \ .
\end{equation}

Letting $x = S$ and $y = ST$ as in our discussion of $\Gamma$,
we see that $2\times2$ matrices for $x$ and $y$ are:
\begin{equation}\label{xy2}
x=\begin{pmatrix}0 & -1\\
1 & 0
\end{pmatrix}\quad,\quad y=\begin{pmatrix}0 & -1\\
1 & \lambda_{n}
\end{pmatrix} \ .
\end{equation}

For $H_{n}$, we clearly have $\left(ST\right)^{n}=I$, thereby constraining $\lambda_{n}$ for a given $n$ \cite{HeckeDessins}. 
Diagonalising $y$ to compute $y^n$ explicitly places a constraint which allows for a solution of $\lambda_n$; we find  the following general expression, as well as some important values for small $n$:
\begin{equation}
\lambda_n =2\cos\left(\pi/n\right) \ , \qquad
\begin{tabular}{|c||c|c|c|c|}
\hline 
$n$ & 3 & 4 & 5 & 6\tabularnewline
\hline 
$\lambda_{n}$ & 1 & $\sqrt{2}$ & $\frac{1+\sqrt{5}}{2}$ & $\sqrt{3}$ \tabularnewline
\hline 
\end{tabular}
\end{equation}
In particular, the $\lambda_n$ are algebraic numbers.

\subsubsection{Congruence Subgroups of Hecke Groups}

By way of extension of the above discussion of subgroups of the modular group $\Gamma$, it is
useful to consider congruence subgroups of Hecke groups. 
The Hecke groups are discrete subgroups
of $\mathrm{PSL}\left(2,\mathbb{R}\right)$; in fact, the matrix entries are in $\mathbb{Z}\left[\lambda_{n}\right]$, the extension of the ring of integers by the algebraic number $\lambda_{n}$. Note that:
\begin{equation}
 H_{n}\subset \mathrm{PSL}\left(2,\mathbb{Z}\left[\lambda_{n}\right]\right) \ . 
\end{equation} 
 However, unlike the special case of the modular group, this inclusion is strict. With this point in mind,
 we can define the congruence subgroups of Hecke groups in the following way
\cite{HeckeDessins}. Let $I$ be an ideal of $\mathbb{Z}\left[\lambda_{n}\right]$. We then define:
\begin{equation}
\mathrm{PSL}\left(2,\mathbb{Z}\left[\lambda_{n}\right],I\right)=\left\{ \begin{pmatrix}a & b\\
c & d
\end{pmatrix}\in\mathrm{PSL}\left(2,\mathbb{Z}\left[\lambda_{n}\right]\right)\;;\;a-1,b,c,d-1\in I\right\} \ . 
\end{equation}

By analogy, we also define:
\begin{align}
\mathrm{PSL^{1}}\left(2,\mathbb{Z}\left[\lambda_{n}\right],I\right)
&=\left\{ \begin{pmatrix}a & b\\
c & d
\end{pmatrix}\in\mathrm{PSL}\left(2,\mathbb{Z}\left[\lambda_{n}\right]\right)\;;\;a-1,c,d-1\in I\right\} 
\\
\mathrm{PSL^{0}}\left(2,\mathbb{Z}\left[\lambda_{n}\right],I\right)
&=\left\{ \begin{pmatrix}a & b\\
c & d
\end{pmatrix}\in\mathrm{PSL}\left(2,\mathbb{Z}\left[\lambda_{n}\right]\right)\;;\;c\in I\right\} 
\ .
\end{align}

Then we can define the congruence subgroups $H_{n}\left(m\right)$, $H_{n}^{1}\left(m\right)$
and $H_{n}^{0}\left(m\right)$ of $H_{n}$ as follows:
\begin{align}
H_{n}\left(I\right)&=\mathrm{PSL}\left(2,\mathbb{Z}\left[\lambda_{n}\right],I\right)\cap H_{n} \ ;
\\
H_{n}^{1}\left(I\right)&=\mathrm{PSL^{1}}\left(2,\mathbb{Z}\left[\lambda_{n}\right],I\right)\cap H_{n} \ ;
\\
H_{n}^{0}\left(I\right)&=\mathrm{PSL^{0}}\left(2,\mathbb{Z}\left[\lambda_{n}\right],I\right)\cap H_{n} \ .
\end{align}
We have:
\begin{equation}
H_{n}\left(I\right)
\subseteq H_{n}^{1}\left(I\right)
\subseteq H_{n}^{0}\left(I\right)
\subseteq H_{n} \ .
\end{equation}
 
By analogy with our discussion of the modular group, we define congruence subgroups 
of level $m$ of the Hecke group $H_{n}$ as subgroups of $H_{n}$ containing
$H_{n}\left(m\right)$ but not any $H_{n}\left(p\right)$ for $p<m$ \cite{HeckeSubs}.
With these details regarding Hecke groups and their subgroups in hand, we can now
consider their connections to clean \emph{dessins d'enfants}.

\subsection{Dessins d'Enfants and Belyi Maps}
A \emph{dessin d'enfant} in the sense of Grothendieck is an ordered pair
$\left\langle X,\mathcal{D}\right\rangle $, where $X$
is an oriented compact topological surface and $\mathcal{D}\subset X$ is a finite graph\textbf{
}satisfying the following conditions \cite{DessinBook}:
\begin{enumerate}
\item $\mathcal{D}$ is connected.
\item $\mathcal{D}$ is \emph{bipartite}, i.e.~consists of only black and
white nodes, such that vertices connected by an edge have different
colours.
\item $X\setminus\mathcal{D}$ is the union of finitely many topological
discs, which we call the \emph{faces}.
\end{enumerate}

\begin{figure}[t!!!]

\begin{center}
\begin{minipage}[t]{0.50\textwidth}%
\begin{center}
\includegraphics[scale=0.2]{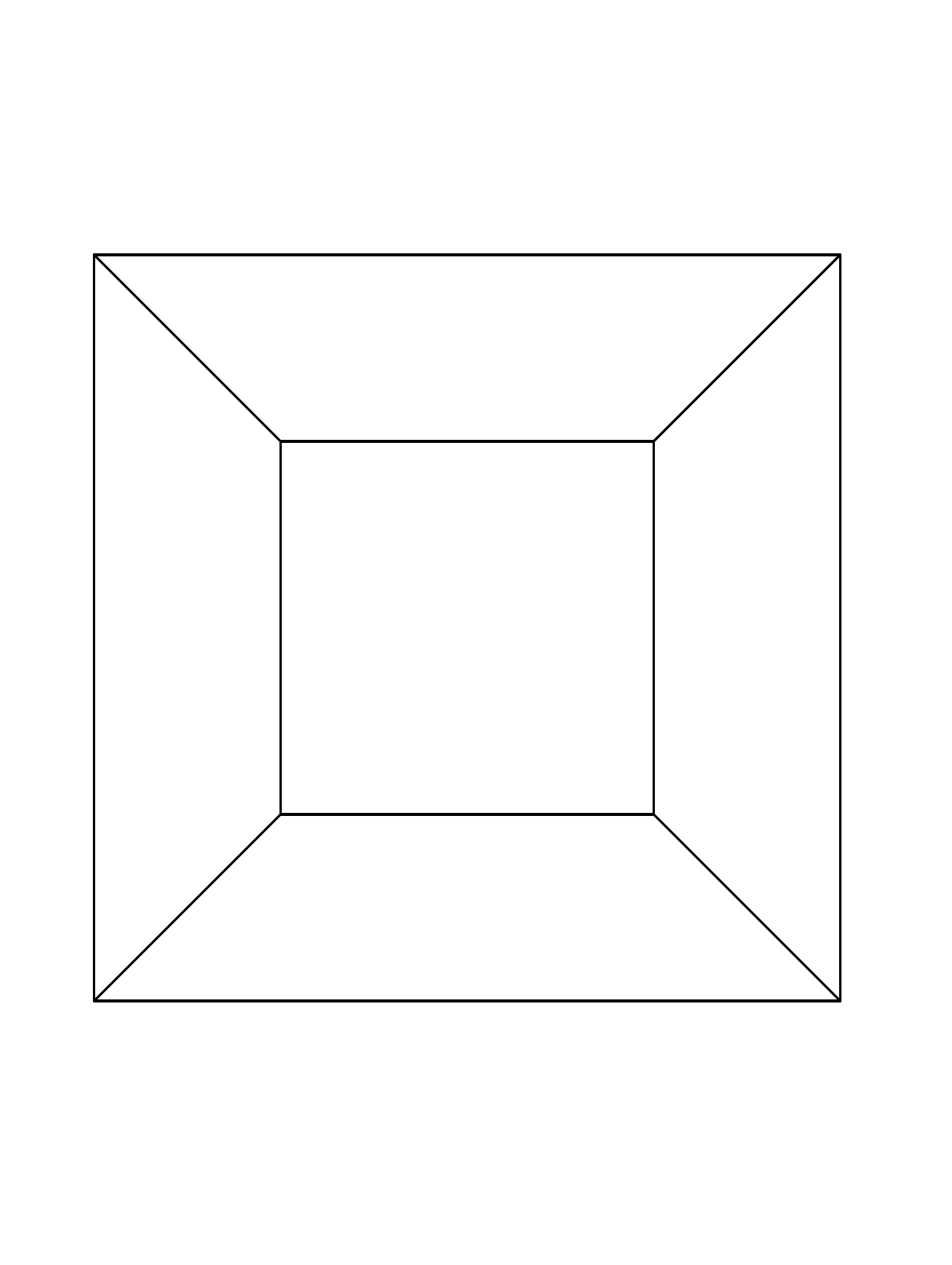}
\par\end{center}
\begin{center}
\small{(a): The planar graph for the cube.}
\par\end{center}%
\end{minipage}%
\begin{minipage}[t]{0.50\textwidth}%
\begin{center}
\includegraphics[scale=0.2]{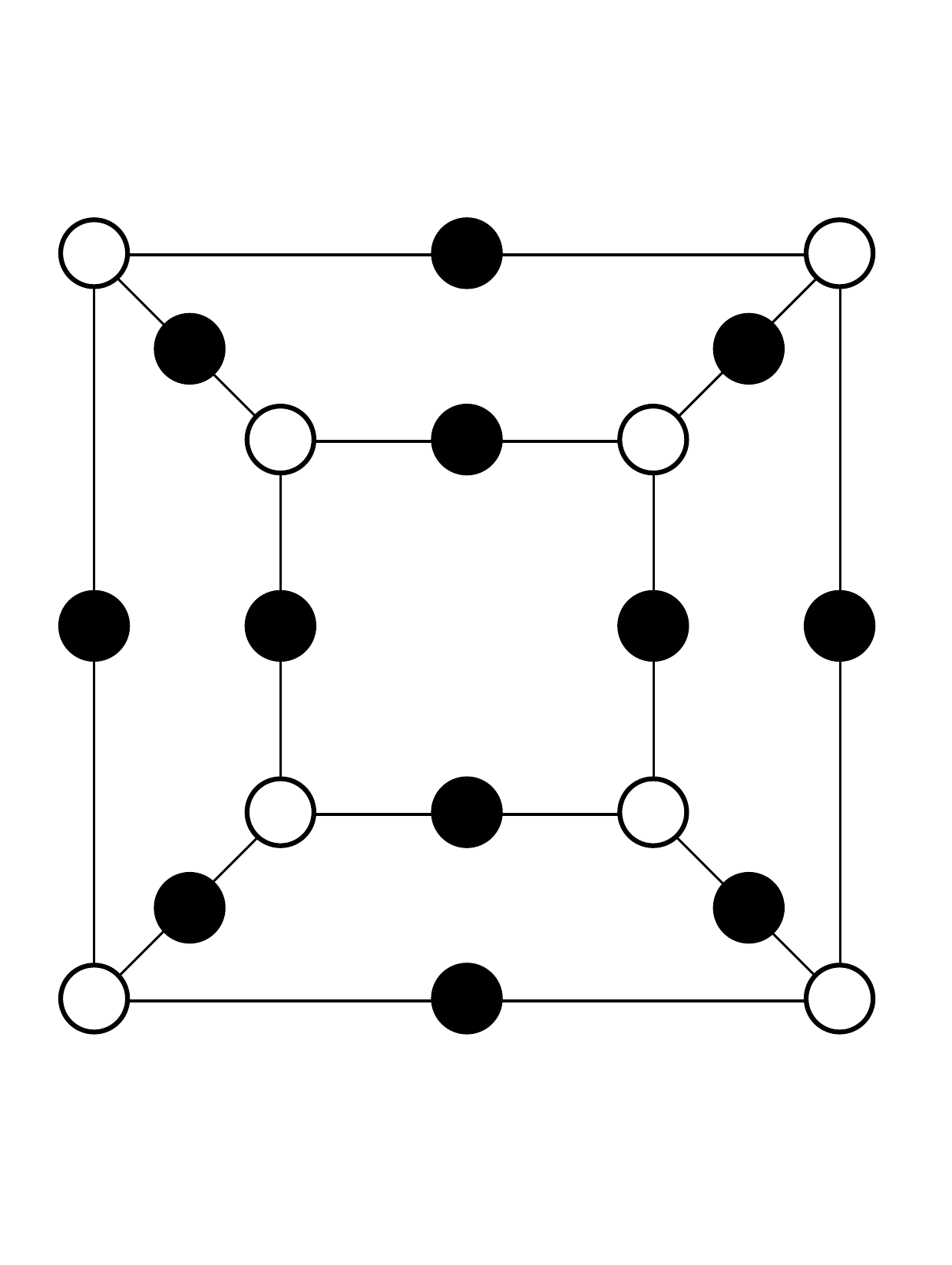}
\par\end{center}
\begin{center}
\small{(b): The corresponding clean \emph{dessin}.}
\par\end{center}%
\end{minipage}%
\end{center}

\sf{\caption{Interpreting the planar graph of the cube as a clean \emph{dessin}.}
\label{f:dessineg}}
\end{figure}

We can interpret any polytope as a \emph{dessin} by inserting a black node into every edge,
and colouring all vertices white.
This process of inserting into each edge a bivalent node of a certain
colour is standard in the study of \emph{dessins d'enfants} and gives rise to so-called \emph{clean} \emph{dessins},
i.e.~those for which all the nodes of one of the two possible colours have valency two.
An example of this procedure for the cube is shown in Figure \ref{f:dessineg}.

There is a one-to-one correspondence
between \emph{dessins d'enfants} and \emph{Belyi maps} \cite{DessinBook, JS1, JS2}.
A Belyi map is a holomorphic map to $\mathbb{P}^{1}$ ramified at
only $\left\{ 0,1,\infty\right\} $, i.e.~for which the only points
$\tilde{x}$ where $\frac{\mathrm{d}}{\mathrm{d}x}\beta\left(x\right)|_{\tilde{x}}=0$
are such that $\beta\left(\tilde{x}\right)\in\left\{ 0,1,\infty\right\} $.
We can associate a Belyi map $\beta\left(x\right)$ to a \emph{dessin} via
its \emph{ramification indices}:~the order of vanishing of the Taylor
series for $\beta\left(x\right)$ at $\tilde{x}$ is the ramification
index $r_{\beta\left(\tilde{x}\right)\in\left\{ 0,1,\infty\right\} }\left(i\right)$
at that $i$th ramification point \cite{1,YMR}. To draw the \emph{dessin}
from the map, we mark one white node for the $i$th pre-image of 0,
with $r_{0}\left(i\right)$ edges emanating therefrom; similarly,
we mark one black node for the $j$th pre-image of 1, with $r_{1}\left(j\right)$
edges. We connect the nodes with the edges, joining only black with
white, such that each face is a polygon with $2r_{\infty}\left(k\right)$
sides \cite{YMR}. 
The converse direction (from \emph{dessins} to Belyi maps) is detailed in e.g.~\cite{DessinBook}, and involves using the \emph{dessin} to construct a so-called \emph{triangle decomposition} of the surface $X$ in which it is embedded, before constructing a unique Belyi function from that decomposition.
For more information on \emph{dessins} and Belyi maps, the reader is referred to \cite{DessinBook, LMS1, LMS2, LMS3}.

\subsection{Belyi Maps from Hecke Subgroups}

Let $G$ be a torsion-free subgroup of the Hecke group $H_n$ of finite index (i.e.~a subgroup that contains no element of finite order other than the identity). Then, the compact surface $X_G$ is obtained from the surface $\mathcal{H}/G$ (where $\mathcal{H}$ is the upper half plane) by adding finitely many points on the boundary real line of $\mathcal{H}$, called \emph{cusps} (in mathematical language:~one for each $G$-orbit of boundary points at which the stabiliser is non-trivial). For example, choosing $n=3$ and letting $G$ be the full modular group $\Gamma$, there is a single cusp, equivalent to any rational number under $G$, namely, the point at infinity -- this is the one-point compactification of the complex plane into the Riemann sphere.
Then (returning to the general case), after carrying out this compactification, the holomorphic projection $\pi:\mathcal{H}/G \rightarrow \mathcal{H}/H_n \cong \mathbb{C}$ extends to a Belyi map $\beta : X_G \rightarrow \mathbb{P}^1$, with each cusp of $X_G$ mapped by $\beta$ to the single cusp $\left\{\infty\right\}$ of $H_n$ which compactifies the plane $\mathcal{H}/H_n$ to $\mathbb{C}\cup \infty = \mathbb{P}^1$ \cite{Harvey}. In turn, this Belyi map will have a unique associated \emph{dessin}, as discussed in the previous section.


\subsection{Schreier Coset Graphs}\label{s:schreier}

As mentioned, 
the Hecke group $H_{n}$ has presentation $\left\langle x,y|x^{2}=y^{n}=I\right\rangle$,
and is thus the free product of cyclic groups $C_{2}=\left\langle x|x^{2}=I\right\rangle$
and $C_{n}=\left\langle y|y^{n}=I\right\rangle$.
Given the free product structure of $H_{n}$, we see that its Cayley graph is an infinite free $n$-valent tree, but with each node replaced by an oriented $n$-gon.
Now, for any \emph{finite index} subgroup $G$ of $H_{n}$, we can quotient the Cayley graph to arrive at a finite graph by associating nodes to right cosets and edges between cosets which are related by action of a group element.
In other words, this graph encodes the permutation representation of $H_{n}$ acting on the right cosets of $G$. 
This is called a {\em Schreier coset graph}, sometimes also referred to in the literature as a {\em Schreier-Cayley coset graph}, or simply a {\em coset graph}.



It is useful consider the permutations induced by the respective actions of
$x$ and $y$ on the cosets of each Hecke subgroup; denote these by $\sigma_{0}$ and $\sigma_{1}$, respectively. We can find a third permutation $\sigma_{\infty}$ by imposing
the following condition, thereby constructing a \emph{permutation
triple} \cite{JS1, JS2}:
\begin{equation}\label{triple}
\sigma_{0}\cdot\sigma_{1}\cdot\sigma_{\infty}=1 \ .
\end{equation}

The permutations $\sigma_{0}$, $\sigma_{1}$ and $\sigma_{\infty}$
give the \emph{permutation representations} of $H_{n}$ on the right cosets
of each subgroup in question. As elements of the symmetric group,
$\sigma_{0}$ and $\sigma_{1}$ can be easily computed from
the Schreier coset graphs by following the procedure
elaborated in \cite{3}, i.e.~by noting that the doubly directed edges
represent an element $x$ of order 2, while the positively oriented
triangles represent an element $y$ of order $n$. Since the graphs are
connected, the group generated by $x$ and $y$ is transitive on the
vertices. Clearly, $\sigma_{0}$ and $\sigma_{1}$ tell us which vertex of the coset graph is sent to which, i.e.~which coset of the Hecke subgroup in question is sent to which by the action of $H_{n}$ on the right cosets of this subgroup.

Now, as detailed in\cite{dessincoset}, there is a direct connection between the Schreier coset graphs and the \emph{dessins d'enfants} for each class of Hecke subgroups:~{the \emph{dessins d'enfants} for a certain conjugacy class of subgroups of a Hecke group $H_{n}$
can be constructed from the Schreier coset graphs by replacing each positively oriented $n$-gon with a white node, and inserting a black node into every edge. Conversely, the Schreier coset graphs can be constructed from the \emph{dessins} by replacing each white node with a positively oriented $n$-gon, and removing the black node from every edge.} Crucially, the \emph{dessin} constructed in this way is the {same} \emph{dessin} as that associated to the Belyi map for the surface $X_G$ as discussed above, where $G$ is the conjugacy class of Hecke subgroups under consideration \cite{dessincoset, Harvey}. 


This last result means that our chain of connections between Hecke groups and \emph{dessins d'enfants} has come full circle:~for a given conjugacy class of subgroups $G$ of a Hecke group $H_n$, there is an associated (compactified) surface $X_G$; in each case there exists a Belyi map $\beta : X_G \rightarrow \mathbb{P}^1$. Every such Belyi map has a unique associated \emph{dessin}, which in turn can be translated into a Schreier coset graph for a certain conjugacy class of Hecke subgroups -- precisely the $G$ from which we began!

\section{Hecke Subgroups and Archimedean Solids}

Having reviewed the theory connecting \emph{dessins d'enfants} with Hecke subgroups, we can now apply this theory to the particular case of the {Archimedean solids}.
In this section, we first give a precise definition of these geometrical entities.
We then identify the conjugacy class of Hecke subgroups corresponding to each Archimedean solid, by giving the permutation representations
$\sigma_{0}$ and $\sigma_{1}$ for each such class. Finally, we discuss some interesting aspects of these results.

\subsection{Platonic and Archimedean Solids}

The Platonic solids are the regular, convex polyhedra; they are the tetrahedron, cube, octahedron, dodecahedron and icosahedron.
 In order to introduce the wider class
of \emph{Archimedean solids}, consider planar graphs without loops, and with
vertices of degree $k>2$. Following \cite{Zvon}, let us call the list of numbers
$\left(f_{1},f_{2},\ldots,f_{k}\right)$, where the $f_{i}$ are the 
number of edges of the adjacent faces taken in 
the counter-clockwise direction around the vertex, the \emph{type} of that vertex. Two such lists
are \emph{equivalent} if one can be obtained from the other by (a) making a cyclic shift and (b)
inverting the order of the $f_{i}$. A solid is called {\em Archimedean} if the types of all its vertices are equivalent \cite{Zvon}.\footnote{It is worth listing some classical sources on the Archimedean solids. Theses objects were first discussed by Archimedes, in a now lost work to which Pappus refers \cite{G}. In the 15\textsuperscript{th} century, the solids were rediscovered by Kepler \cite{F}. Classical geometers who discuss the Archimedean solids include Sommerville \cite{Sommerville} and Miller \cite{G}.}

We emphasise a subtlety here. One informal way of defining the Archimedean
solids is as the semi-regular convex polyhedra composed of two or more types of regular polygons meeting in
identical vertices, with no requirement that faces be equivalent.
\emph{Identical vertices} is usually taken to mean that for any two vertices, there must be an isometry of the entire solid that takes one vertex to the other. Sometimes, however, it is only required that the faces that meet at one vertex are related isometrically to the faces that meet at the other. On the former definition, the so-called \emph{pseudorhombicuboctahedron}, also known as the \emph{elongated square gyrobicupola}, is not considered an Archimedean solid; on the latter it is. The formal definition
of the Archimedean solids above corresponds to the latter definition here; it is this latter definition which which we shall use.

\begin{table}[t]
\resizebox{\textwidth}{!}{%
\scriptsize
\begin{tabular}{|c||c|c|c|c|c|c|c|c|}
\hline 
& Name & $V$ & Vertex type & $F$ & $E$ & Sym & Hecke Group \tabularnewline
\hline 
\hline 
I & Tetrahedron & 4 & $\left(3,3,3\right)$ & 4 & 6 & $T_{d}$ & $H_{3}$ \tabularnewline
\hline
I & Cube & 8 & $\left(4,4,4\right)$ & 6 & 12 & $O_{h}$ & $H_{3}$ \tabularnewline
\hline
I & Octahedron & 6 & $\left(3,3,3,3\right)$ & 8 & 12 & $O_{h}$ & $H_{4}$ \tabularnewline
\hline
I & Dodecahedron & 20 & $\left(5,5,5\right)$ & 12 & 30 & $I_{h}$ & $H_{3}$ \tabularnewline
\hline
I & Icosahedron & 12 & $\left(3,3,3,3,3\right)$ & 20 & 30 & $I_{h}$ & $H_{5}$ \tabularnewline
\hline
II & $n$-prism & $2n$ & $\left(4,4,n\right)$ & $n+2$ & $3n$ & $D_{nh}$ & $H_{3}$ \tabularnewline
\hline 
II & $n$-antiprism & $2n$ & $\left(3,3,3,n\right)$ & $2n+2$ & $4n$ & $D_{nd}$ & $H_{4}$ \tabularnewline
\hline 
III & Truncated tetrahedron & 12 & $\left(3,6,6\right)$ & 8 & 18 & $T_{d}$ & $H_{3}$ \tabularnewline
\hline 
III & Truncated cube & 24 & $\left(3,8,8\right)$ & 14 & 36 & $O_{h}$ & $H_{3}$ \tabularnewline
\hline 
III & Truncated octahedron & 24 & $\left(4,6,6\right)$ & 14 & 36 & $O_{h}$ & $H_{3}$ \tabularnewline
\hline 
III & Truncated isocahedron & 60 & $\left(5,6,6\right)$ & 32 & 90 &  $I_{h}$ & $H_{3}$ \tabularnewline
\hline 
III & Truncated dodecahedron & 60 & $\left(3,10,10\right)$ & 32 & 90 &  $I_{h}$ & $H_{3}$ \tabularnewline
\hline 
III & Truncated cuboctahedron & 48 & $\left(4,6,8\right)$ & 26 & 72 & $O_{h}$ & $H_{3}$ \tabularnewline
\hline 
III & Truncated icosidodecahedron & 120 & $\left(4,6,10\right)$ & 62 & 180 &  $I_{h}$ & $H_{3}$ \tabularnewline
\hline 
III & Cuboctahedron & 12 & $\left(3,4,3,4\right)$ & 14 & 24 & $O_{h}$ & $H_{4}$ \tabularnewline
\hline 
III & Icosidodecahedron & 30 & $\left(3,5,3,5\right)$ & 32 & 60 &  $I_{h}$ & $H_{4}$ \tabularnewline
\hline 
III & Rhombicuboctahedron & 24 & $\left(3,4,4,4\right)$ & 26 & 48 & $O_{h}$ & $H_{4}$ \tabularnewline
\hline 
III & Rhombicosidodecahedron & 60 & $\left(3,4,5,4\right)$ & 62 & 120 &  $I_{h}$ & $H_{4}$ \tabularnewline
\hline 
III & Pseudorhombicuboctahedron & 24 & $\left(3,4,4,4\right)$ & 26 & 48 & $D_{4d}$ & $H_{4}$ \tabularnewline
\hline 
III & Snub cube & 24 & $\left(3,3,3,3,4\right)$ & 38 & 60 & $O$ & $H_{5}$ \tabularnewline
\hline 
III & Snub dodecahedron & 60 & $\left(3,3,3,3,5\right)$ & 92 & 150 & $I$ & $H_{5}$ \tabularnewline
\hline 
\end{tabular}}
\sf{\caption{The Archimedean solids.
Type I are the five Platonic solids.
Type II are the two infinite families:~the prisms and anti-prisms.
Type III are the fourteen exceptional solids.
We record the number of vertices $V$, faces $F$ and edges $E$ in each case.
The \emph{vertex type} is the number of edges of the adjacent faces to the vertex, taken counter-clockwise; the fact that all vertices have the same vertex type is the definition of \emph{Archimedean}.
`Sym' denotes the symmetry group of the solid, given in Sch\"onflies notation.
We also indicate the associated Hecke group, determined by the valency of the vertices of each solid.
\label{ArchTable}
}}
\end{table}

The solids which satisfy the condition for being \emph{Archimedean} are:
\begin{itemize}
\item[I] The five Platonic solids; 
\item[II] Two infinite series (the prisms and anti-prisms);
\item[III] Fourteen exceptional solutions.\footnote{An interesting aside:~It is known that the five Platonic solids correspond to three symmetry groups, which in turn are related to the exceptional Lie algebras $E_{6,7,8}$, and furthermore to Arnold's simple surface singularities with no deformations \cite{BPR}.
It is curious that the next order generalisation of the simple singularities, with exactly one deformation modulus, has fourteen exceptional cases. We conjecture, therefore, some connection to the exceptional Archimedean
solids.
}
\end{itemize}
All these solids are listed in table \eqref{ArchTable}.
For completeness, we also give the symmetry groups of each of these solids in the column labelled `Sym', written in {\em Sch\"onflies notation}. Here, we recognise the standard tetrahedral, octahedral and icosahedral groups:
\begin{equation}
T \simeq A_4 \ , \quad
O \simeq S_4 \ , \quad
I \simeq  A_5 \ ,
\end{equation}
as well as the dihedral group $D_n$.
The subscripts $d$ and $h$ denote extra symmetries about a horizontal ($h$) or diagonal ($d$) plane.

\subsection{Hecke Subgroups for the Archimedean Solids}

As discussed, we can interpret all the Archimedean solids
as clean \emph{dessins d'enfants} by inserting a black node into every edge and colouring every vertex white.
By the correspondences detailed in the previous section, we should then be able to associate
a class of Hecke subgroups to every Archimedean solid.
To achieve this, we begin by converting the \emph{dessins} for these solids to Schreier coset graphs, reading off the associated  permutations $\sigma_{0}$ and $\sigma_{1}$
for the conjugacy class of subgroups of the relevant Hecke group in each case; these results are tabulated in appendix A. Next, we find a representative of the conjugacy class of subgroups
of the relevant Hecke group in each case by inputting the permutations $\sigma_{0}$ and $\sigma_{1}$ into the \textsf{GAP} \cite{GAP} algorithm presented in appendix B. This yields a representative of the conjugacy classes of subgroups in each case.

To illustrate this procedure, consider the specific example of the octahedron. Inputting the permutations for the octahedron into the algorithm of appendix B yields the following generators for a representative of the
conjugacy class of subgroups of $H_{4}$ associated to this solid:

\begin{align}
& \left\{
\notag
\left(
\begin{array}{cc}
 -1 & 0 \\
 3 \sqrt{2} & -1 \\
\end{array}
\right),\left(
\begin{array}{cc}
 1 & -3 \sqrt{2} \\
 0 & 1 \\
\end{array}
\right),\left(
\begin{array}{cc}
 -7 & -3 \sqrt{2} \\
 6 \sqrt{2} & 5 \\
\end{array}
\right),\left(
\begin{array}{cc}
 -7 & -6 \sqrt{2} \\
 3 \sqrt{2} & 5 \\
\end{array}
\right),
\right.
\\
& \left.
\left(
\begin{array}{cc}
 5 & -6 \sqrt{2} \\
 3 \sqrt{2} & -7 \\
\end{array}
\right),\left(
\begin{array}{cc}
 5 & -3 \sqrt{2} \\
 6 \sqrt{2} & -7 \\
\end{array}
\right),\left(
\begin{array}{cc}
 11 & -15 \sqrt{2} \\
 -15 \sqrt{2} & 41 \\
\end{array}
\right)
\right\} \ .
\end{align}

Note that comparison with our results concerning congruence subgroups of Hecke groups in section \eqref{s:hecke} reveals that the octahedron corresponds to the principal congruence subgroup $H_{4}\left(3\right)$.
Repeating this procedure for all Archimedean solids, we find a number of results worthy of comment:

\begin{itemize}
\item While the octahedron corresponds to the conjugacy class of Hecke subgroups $H_{4}\left(3\right)$, no other Archimedean solid is readily associated with congruence Hecke subgroups in an analogous manner.

\item From section \eqref{s:hecke}, the matrices of $H_{4}$ are of the following two types:

\begin{equation}
\begin{pmatrix}a & b\sqrt{2}\\
c\sqrt{2} & d
\end{pmatrix}\quad,\quad\begin{pmatrix}a\sqrt{2} & b\\
c & d\sqrt{2}
\end{pmatrix}
\end{equation}

The elements of the first type form a subgroup of index 2 in $H_{4}$; all the subgroups of $H_{4}$
corresponding to the Archimedean solids turn out to be subgroups of this subgroup.

\item We have seen that, interpreted as \emph{dessins}, the tetrahedron, cube and 
dodecahedron correspond to the conjugacy classes of subgroups
$\Gamma\left(3\right)$, $\Gamma\left(4\right)$ and $\Gamma\left(5\right)$, respectively \cite{YMR}. In addition to these three arising as the \emph{dessins} corresponding to certain congruence subgroups of the modular group, the \emph{dessins} for the 33
conjugacy classes of genus zero,
torsion-free congruence subgroups (discussed in detail in \cite{YMR}) reveal that \emph{other} Archimedean solids also correspond to
certain conjugacy classes of congruence subgroups of the modular group. Specifically, the truncated tetrahedron is associated with the modular subgroup
$\Gamma_{0}\left(2\right)\cap\Gamma\left(3\right)$, while the 8-prism is associated with the subgroup
$\Gamma\left(8;2,1,2\right)$.
Given this, the question arises as to whether the other trivalent Archimedean solids can also be 
associated with congruence subgroups of the modular group. This can be checked in \textsf{Sage} \cite{sage} using the permutations $\sigma_0$ and $\sigma_1$ as input; doing so, we find that the answer to this question is  in general \emph{negative}. As a more specific result, we find that the truncated tetrahedron is the \emph{only} trivalent, exceptional Archimedean solid which corresponds to a conjugacy class of congruence subgroups of the modular group.

\end{itemize}

\section{Conclusions}

In this paper, we have reviewed the connections between \emph{dessins d'enfants} and Hecke subgroups.
We have applied this theory to the case of the Archimedean solids, showing how, through interpretation as clean \emph{dessins}, each of these geometric objects can be associated with a conjugacy class of subgroups of Hecke groups. In addition to opening up one more area of research into the Archimedean solids, this work is useful from the point of view of mathematical physics, as it can help to shed new light on cases where \emph{dessins} naturally arise in physical contexts. To take two examples:

\begin{enumerate}

\item Certain $\mathcal{N}=2$ supersymmetric gauge theories in four dimensions \cite{Gaiotto,Dummies} naturally give rise to trivalent \emph{dessins} \cite{HR2015} (including {Archimedean} \emph{dessins}); with the theory of this paper in mind, this suggests an otherwise overlooked connection between these gauge theories and modularity. 

\item In \cite{YMR,1,He:2013lha}, it is shown that certain Calabi-Yaus give rise to clean, trivalent \emph{dessins}, including many Archimedean \emph{dessins}. Again, the theory of this paper therefore suggests a connection between the theory of Hecke groups and these Calabi-Yaus.

\end{enumerate}

The general lessson is the following:~given the growing ubiquity of \emph{dessins} in mathematical physics, it is important to be clear and explicit on the connections of these objects to other areas of mathematics. Hence, it is important to illustrate the theory underlying these connections. It is particularly useful to lay out this theory explicitly in special cases such as that of the Archimedean solids, which we expect (on both inductive grounds and grounds of simplicity) to arise with greatest frequency in physical contexts.

\section*{Acknowledgements}

Y-H.H.~thanks the Science and Technology Facilities Council, UK, for an Advanced Fellowship and for STFC grant ST/J00037X/1; the Chinese Ministry of Education, for a Chang-Jiang Chair Professorship at NanKai University; the city of Tian-Jin for a Qian-Ren Scholarship; the US NSF for grant CCF-1048082; as well as City University, London; the Department of Theoretical Physics, Oxford; and Merton College, Oxford, for their enduring support. J.R.~is supported by an AHRC studentship at 
the University of Oxford, and is also indebted to Merton College for their support.


\begin{appendix}

\section{Permutations for the Archimedean Solids}

In this appendix, we present the permutations $\sigma_0$ and $\sigma_1$ corresponding to the conjugacy classes of Hecke subgroups for each of the Archimedean solids. Although it is in principle easy to read these permutations off from the Schreier coset graphs in the manner detailed in section \eqref{s:schreier}, the procedure is extremely time-consuming, so it is worth presenting the results in full here. The method for computing explicit generators for a representative of each class of subgroups is detailed in the following appendix.

As a technical point, note that we have just presented the permutations $\sigma_0$ in each case. The permutations $\sigma_1$ can be written down by populating $2E/n$ $n$-tuples with numerals in ascending order from $1$ to $2E$, where $E$ is the number of edges of the solid, and $H_n$ is the associated Hecke group as given in table \eqref{ArchTable}. So, for example, the permutations $\sigma_1$ for the octahedron are: (1, 2, 3, 4) (5, 6, 7, 8) (9, 10, 11, 12) (13, 14, 15, 16) (17, 18, 19, 20) (21, 22, 23, 24).

\subsection{Platonic Solids}

Interpreted as clean \emph{dessins d'enfants}, the tetrahedron, cube, and 
dodecahedron correspond to the conjugacy classes of principal congruence subgroups
$\Gamma\left(3\right)$, $\Gamma\left(4\right)$ and $\Gamma\left(5\right)$, respectively \cite{YMR}. In addition, from table \eqref{ArchTable}
we can see that the remaining two Platonic solids -- the octahedron and the icosahedron -- correspond to conjugacy classes of subgroups of $H_{4}$ and $H_{5}$, respectively.\footnote{Of course, one should note that the octahedron and icosahedron \emph{could} also be defined as subgroups of the modular group $\Gamma$, as their duals are trivalent.} The permutations $\sigma_0$ for these classes are:

\paragraph{Octahedron:} (1, 5) (6, 11) (4, 12) (2, 13) (3, 23) (8, 14) (7, 18) (10, 19) (9, 22) (16, 24) (15, 17) (20, 21)


\paragraph{Icosahedron:} (1, 6) (10, 11) (5, 15) (2, 21) (3, 16) (4, 43) (7, 22) (8, 27) (9, 33) (12, 34) (13, 37) (14, 42) (17, 25) (23, 26) (28, 32) (35, 36) (38, 41) (20, 44) (18, 46) (24, 47) (30, 48) (29, 52) (31, 53) (40, 54) (39, 57) (45, 58) (19, 59) (50, 60) (55, 56) (49, 51)
 

\subsection{Prisms and Antiprisms}
\label{sec:prisms}

The prisms and antiprisms form infinite series of Archimedean solids. Since all prisms are
trivalent, all correspond to a conjugacy class of subgroups of $\Gamma$; since all antiprisms
are 4-valent, all correspond to a conjugacy class of subgroups of $H_{4}$. We
{can} construct general expressions for the permutations $\sigma_{0}$ and $\sigma_{1}$
for the prisms and antiprisms. These expressions can be used to find the permutations 
$\sigma_{0}$ and $\sigma_{1}$ for any particular (anti)prism of interest.

\paragraph{Prisms:} Call the $n$-gon faced prism the $n$-prism. The $n$-prism has the following
permutations $\sigma_{0}$ and $\sigma_{1}$:
\begin{align}
\notag \sigma_{0}:\ & \left(3n-2,3\right) \left(3n-1,6n-2\right)\left(6n-1,3n+3\right)\\
\notag & \cdot\prod\limits_{i=0}^{n-2}\left(3i+1,3i+6\right)\left(3i+2,3n+3i+1\right)\left(3n+3i+2,3n+3i+6\right)\\
\label{prism}
\sigma_{1}:\ & \prod\limits_{i=0}^{n+2}\left(3i+1,3i+2,3i+3\right) \ .
\end{align}

\paragraph{Antiprisms:} Call the $n$-gon faced anti-prism the $n$-antiprism. The $n$-antiprism has
 the following permutations $\sigma_{0}$ and $\sigma_{1}$:

\begin{align}
\notag \sigma_{0}:\ & \left(4n-3,4\right) \left(4n-2,8n-2\right) \left(3,8n-1\right) \left(4n+1,8n\right)\\
\notag & \cdot\prod\limits_{i=0}^{n-2}\left(4i+1,4i+8\right)\left(4i+2,4n+4i+2\right)\left(4i+7,4n+4i+3\right)\left(4n+4i+4,4n+4i+5\right)\\
\sigma_{1}:\ & \prod\limits_{i=0}^{n+2}\left(4i+1,4i+2,4i+3,4i+4\right) \ .
\end{align}

\subsection{Exceptional Archimedean Solids}

In addition to the Platonic solids and the prisms and antiprisms, there remain fourteen exceptional
Archimedean solids, as given in table \eqref{ArchTable}. Here we give the permutations $\sigma_0$ for each corresponding class of Hecke subgroups.

\paragraph{Truncated tetrahedron:} (1, 32) (3, 4) (2, 7) (5, 9) (6, 24) (8, 10) (11, 13) (12, 18) (15, 16) (17, 19) (14, 28) (29, 31) (33, 34) (30, 36) (26, 35) (20, 25) (23, 27) (21, 22)

\paragraph{Truncated cube:} (3, 4) (2, 7) (5, 9) (60, 61) (63, 55) (57, 58) (20, 22) (21, 27) (24, 25) (11, 13) (12, 18) (15, 16) (39, 45) (43, 42) (38, 40) (30, 31) (29, 35) (32, 34) (66, 72) (65, 68) (67, 70) (48, 49) (47, 53) (50, 52) (1, 59) (62, 64) (69, 54) (51, 6) (41, 46) (8, 10) (56, 23) (36, 71) (33, 37) (17, 44) (14, 19) (26, 28)

\paragraph{Truncated octahedron:} (1, 11) (2, 4) (5, 8) (9, 10) (12, 13) (3, 36) (6, 56) (7, 65) (35, 26) (32, 34) (29, 31) (27, 28) (25, 17) (16, 14) (15, 24) (21, 22) (18, 19) (23, 68) (70, 69) (66, 67) (62, 64) (63, 72) (59, 61) (57, 58) (51, 60) (50, 52) (54, 55) (49, 48) (45, 71) (20, 42) (30, 37) (39, 40) (38, 47) (44, 46) (41, 43) (33, 53)

\paragraph{Truncated icosahedron:} (1, 14) (11, 13) (8, 10) (5, 7) (2, 4) (3, 16) (15, 28) (12, 26) (9, 22) (6, 20) (19, 35) (32, 34) (17, 31) (18, 58) (56, 60) (29, 55) (30, 53) (50, 52) (27, 49) (25, 47) (45, 46) (44, 23) (24, 41) (38, 40) (21, 37) (36, 73) (33, 64) (59, 61) (57, 110) (54, 107) (51, 98) (48, 95) (43, 86) (42, 83) (39, 77) (71, 75) (70, 68) (65, 67) (66, 62) (63, 116) (115, 113) (111, 112) (108, 109) (106, 104) (101, 103) (100, 99) (96, 97) (92, 94) (89, 91) (87, 88) (84, 85) (80, 82) (119, 79) (78, 118) (74, 76) (120, 122) (72, 164) (69, 158) (117, 155) (114, 149) (146, 105) (140, 102) (93, 137) (131, 90) (128, 81) (163, 161)  (159, 160) (156, 157) (152, 154) (150, 151) (147, 148) (143, 145) (141, 142) (138, 139) (134, 136) (132, 133) (129, 130) (125, 127) (123, 124) (121, 165) (162, 166) (153, 179) (144, 176) (135, 173) (126, 170) (167, 169) (171, 172) (174, 175) (177, 178) (168, 180)

\paragraph{Truncated dodecahedron:} (1, 29) (3, 32) (30, 31) (27, 44) (23, 25) (24, 43) (21, 41) (18, 40) (17, 19) (11, 13) (15, 39) (12, 38) (5, 7) (9, 35) (6, 34) (2, 4) (26, 28) (22, 20) (14, 16) (8, 10) (36, 62) (33, 46) (45, 107)  (42, 92) (37, 77) (61, 59) (60, 66) (63, 64) (47, 49) (51, 120) (48, 119) (108, 109) (111, 105) (104, 106) (90, 96) (93, 94) (89, 91) (74, 76) (78, 79) (75, 81) (58, 56) (50, 52) (55, 53) (54, 122) (57, 123) (121, 136) (117, 118) (112, 110) (116, 113) (115, 135) (114, 134) (133, 149) (101, 103) (95, 97) (98, 100) (102, 131) (99, 130) (132, 146) (86, 88) (80, 82) (83, 85) (87, 128) (84, 127) (129, 142) (71, 73) (65, 67) (68, 70) (72, 126) (69, 125) (124, 139) (137, 154) (156, 153) (152, 138) (151, 180) (179, 150) (148, 176) (177, 178) (174, 175) (173, 147) (170, 145) (171, 172) (168, 169) (165, 166) (164, 144) (143, 167) (162, 163) (159, 160) (141, 158) (140, 161) (155, 157)

\paragraph{Truncated cuboctahedron:} (1, 23) (24, 20) (21, 17) (18, 14) (15, 12) (10, 9) (7, 6) (4, 3) (2, 25) (22, 29) (19, 32) (16, 35) (13, 38) (11, 41) (8, 44) (5, 46) (26, 28) (30, 55) (56, 58) (31, 59) (33, 34) (36, 62) (63, 64) (65, 37) (39, 40) (42, 67) (68, 70) (43, 71) (45, 48) (47, 49) (50, 52) (27, 53) (54, 77) (51, 73) (76, 74) (78, 79) (75, 119) (57, 86) (60, 89) (88, 87) (85, 83) (90, 91) (61, 98) (66, 101) (100, 99) (95, 97) (102, 103) (72, 113) (69, 110) (111, 112) (114, 115) (107, 109) (80, 82) (84, 128) (126, 127) (81, 125) (129, 130) (132, 133) (96, 134) (92, 94) (93, 131) (135, 136) (137, 105) (104, 106) (108, 140) (138, 139) (141, 142) (117, 143) (123, 144) (120, 121) (116, 118) (122, 124)

\paragraph{Truncated icosidodecahedron:} (1, 29) (26, 28) (23, 25) (20, 22) (17, 19) (15, 16) (11, 14) (8, 10) (5, 7) (2, 4) (3, 31) (30, 59) (27, 56) (24, 52) (21, 50) (18, 47) (13, 44) (12, 41) (9, 38) (6, 34) (33, 60) (58, 89) (88, 86) (85, 57) (53, 55) (83, 54) (80, 82) (79, 51) (48, 49) (46, 77) (75, 76) (45, 74) (42, 43) (40, 71) (68, 70) (39, 67) (35, 37) (36, 65) (64, 62) (61, 32) (90, 119) (87, 116) (84, 113) (81, 109) (78, 107) (73, 104) (72, 101) (69, 98) (66, 95) (63, 91) (93, 121) (123, 179) (178, 176) (175, 174) (173, 120) (118, 117) (115, 170) (169, 167) (166, 164) (163, 162) (161, 114) (112, 110) (111, 158) (155, 157) (152, 154) (149, 151) (108, 148) (105, 106) (103, 146)  (143, 145) (140, 142) (138, 139) (102, 137) (99, 100) (97, 134) (131, 133) (128, 130) (125, 127) (96, 124) (92, 94) (122, 181) (180, 269) (177, 260) (257, 172) (171, 254) (168, 251) (165, 242) (160, 239) (159, 236) (156, 233) (153, 224) (150, 221) (147, 218) (144, 215) (141, 206) (136, 203) (135, 200) (132, 197) (129, 188) (126, 185) (184, 182) (183, 270) (268, 266) (265, 262) (264, 261) (259, 258) (256, 255) (252, 253) (248, 250) (247, 245) (244, 243) (241, 240) (237, 238) (234, 235) (230, 232) (227, 229) (226, 225) (222, 223) (219, 220) (216, 217) (214, 212) (209, 211) (207, 208) (204, 205) (202, 201) (199, 198) (196, 194) (191, 193) (189, 190) (186, 187) (267, 271) (263, 323) (249, 314) (246, 311) (231, 302) (228, 299) (290, 213) (210, 287) (195, 278) (192, 274) (276, 329) (328, 326) (325, 272) (273, 324) (322, 320) (319, 317) (316, 315) (313, 312) (308, 310) (307, 305) (304, 303) (300, 301) (296, 298) (293, 295) (291, 292) (288, 289) (284, 286) (281, 283) (279, 280) (275, 277) (330, 332) (327, 359) (321, 356) (353, 318) (350, 309) (347, 306) (344, 297) (294, 341) (285, 338) (282, 335) (333, 334) (331, 360) (358, 357) (355, 354) (351, 352) (348, 349) (345, 346) (342, 343) (339, 340) (337, 336)

\paragraph{Cuboctahedron:} (1, 10) (11, 14) (15, 8) (4, 5) (2, 21) (3, 17) (20, 6) (7, 31) (30, 16) (27, 13) (22, 9) (12, 26) (25, 39) (28, 42) (43, 29) (32, 46) (47, 19) (18, 33) (24, 34) (23, 38) (35, 37) (40, 41) (44, 45) (36, 48)

\paragraph{Icosidodecahedron:} (1, 9) (10, 17) (18, 26) (27, 34) (35, 4) (3, 37) (2, 5) (6, 12) (11, 13) (14, 20) (19, 21) (22, 25) (28, 29) (30, 33) (36, 40) (38, 41) (44, 79) (39, 78) (8, 45) (7, 49) (46, 52) (16, 54) (55, 57) (15, 58) (24, 63) (64, 65) (23, 66) (32, 70) (31, 74) (71, 73) (75, 77) (80, 110) (76, 109) (42, 48) (47, 82) (81, 43) (50, 53) (56, 86) (51, 85) (60, 93) (59, 62) (61, 94) (67, 69) (72, 102) (101, 68) (89, 96) (95, 100) (90, 99) (98, 107) (103, 108) (97, 104) (105, 112) (111, 116) (106, 115) (84, 113) (83, 117) (114, 120) (88, 118) (87, 92) (91, 119)

\paragraph{Rhombicuboctahedron:} (1, 5) (6, 9) (10, 13) (4, 14) (2, 17) (3, 30) (20, 31) (8, 33) (7, 38) (34, 37) (12, 54) (11, 59) (55, 58) (15, 70) (67, 69) (16, 66) (18, 36) (35, 48) (21, 47) (19, 24) (40, 41) (39, 53) (50, 56) (42, 49) (57, 64) (61, 79) (68, 80) (65, 60) (72, 76) (28, 73) (25, 32) (29, 71) (22, 82) (27, 81) (23, 26) (44, 45) (46, 86) (43, 87) (51, 63) (62, 91) (52, 90) (75, 77) (74, 95) (78, 94) (84, 96) (83, 85) (88, 89) (92, 93)

\paragraph{Rhombicosidodecahedron:} (1, 19) (4, 5) (8, 12) (11, 16) (15, 20) (2, 21) (3, 28) (6, 29) (7, 36) (9, 40) (10, 44) (13, 47) (14, 52) (17, 56) (18, 237) (240, 22) (24, 25) (27, 30) (32, 33) (35, 37) (39, 41) (43, 48) (46, 49) (51, 53) (55, 238) (54, 116) (57, 239) (23, 68) (26, 72) (31, 80) (34, 84) (38, 92) (42, 96) (45, 104) (50, 108) (115, 58) (60, 61) (64, 65) (67, 69) (71, 73) (76, 77) (79, 81) (83, 85) (88, 89) (91, 93) (95, 97) (100, 101) (103, 105) (107, 109) (112, 113) (114, 117) (175, 118) (62, 124) (63, 128) (66, 132) (70, 131) (74, 136) (75, 140) (78, 144) (82, 143) (86, 148) (87, 152) (90, 156) (94, 155) (98, 160) (99, 164) (102, 168) (106, 167) (110, 172) (111, 176) (119, 121) (123, 125) (127, 129) (130, 133) (135, 137) (139, 141) (142, 145) (147, 149) (151, 153) (154, 157) (159, 161) (163, 165) (166, 169) (171, 173) (59, 120) (122, 184) (126, 188) (134, 192) (138, 196) (146, 200) (150, 204) (158, 208)  (162, 212) (170, 213) (174, 177) (183, 185) (187, 189) (191, 193) (195, 197) (199, 201) (203, 205) (207, 209) (211, 214) (178, 216) (180, 181) (182, 224) (186, 223) (190, 228) (194, 227) (198, 232) (202, 231) (206, 236) (210, 235) (215, 218) (179, 217) (222, 225) (226, 229) (230, 233) (234, 219) (220, 221)

\paragraph{Pseudorhombicuboctahedron:} (1, 5) (6, 10) (11, 14) (4, 15) (2, 17) (3, 46) (20, 47) (8, 21) (7, 26) (22, 25) (9, 30) (12, 35) (31, 34)  (13, 38) (16, 43) (39, 42) (18, 24) (23, 53) (50, 56) (19, 49) (28, 58) (59, 61) (32, 62) (27, 29) (33, 65) (66, 70) (40, 71) (36, 37) (44, 45) (48, 80) (76, 79) (41, 75) (52, 77) (51, 81) (84, 96) (78, 95) (55, 82) (83, 88) (60, 85) (54, 57) (86, 64) (63, 68) (67, 91) (87, 90) (89, 93) (73, 94) (72, 74) (69, 92)

\paragraph{Snub cube:} (1, 18) (13, 17) (8, 12) (2, 7) (5, 22) (4, 60) (3, 53) (21, 19) (20, 26) (16, 35) (34, 14) (40, 15) (45, 11) (9, 44) (10, 49) (6, 54) (52, 56) (59, 23) (25, 27) (30, 31) (33, 36) (39, 41) (43, 50) (48, 55) (51, 73) (57, 72) (58, 67) (24, 61) (65, 28) (29, 97) (32, 93) (37, 92) (38, 87) (42, 84) (46, 83) (47, 78) (71, 68) (66, 62) (64, 98) (96, 94) (88, 91) (85, 86) (79, 82) (74, 77) (75, 108) (69, 107) (70, 103) (63, 102) (99, 101) (100, 119) (118, 95) (89, 117) (90, 114) (81, 113) (80, 112) (76, 109) (104, 106) (105, 120) (115, 116) (110, 111)

\paragraph{Snub dodecahedron:} (1, 7) (8, 12) (13, 16) (17, 21) (5, 22) (2, 37) (3, 32) (4, 27) (6, 38) (10, 42) (9, 47) (11, 48) (15, 52) (14, 57) (20, 58) (19, 62) (18, 67) (25, 68) (24, 72) (26, 23) (28, 31) (33, 36) (39, 41) (43, 46) (49, 51) (53, 56) (59, 61) (63, 66) (69, 71) (30, 73) (29, 76) (35, 87) (34, 92) (40, 93) (45, 102) (44, 107) (50, 108) (55, 117) (54, 122) (60, 123) (65, 132) (64, 137) (70, 138) (75, 147) (74, 80) (77, 81) (82, 86) (88, 91) (94, 96) (97, 101) (103, 106) (109, 111) (112, 116) (118, 121) (124, 126) (127, 131) (133, 136) (139, 141) (142, 146) (79, 148) (149, 182) (78, 183) (85, 184) (84, 187) (83, 192) (90, 193) (89, 197) (95, 198) (100, 199) (99, 202) (98, 206) (105, 207) (104, 212) (110, 213) (115, 214) (114, 217) (113, 223) (120, 224) (119, 152) (125, 153) (130, 154) (129, 157) (128, 162) (135, 163) (134, 167) (140, 168) (145, 169) (144, 172) (143, 177) (150, 178) (200, 201) (203, 210) (208, 211) (215, 216) (218, 222) (225, 151) (155, 156) (158, 161) (164, 166) (170, 171) (173, 176) (179, 181) (185, 186) (188, 191) (194, 196) (195, 258) (205, 262) (204, 267) (209, 268) (220, 272) (219, 230) (221, 226) (160, 232) (159, 237) (165, 238) (175, 242) (174, 247) (180, 248) (190, 252) (189, 257) (259, 261) (263, 266) (269, 271) (273, 229) (227, 231) (233, 236) (239, 241) (243, 246)  (249, 251) (253, 256) (260, 293) (265, 294) (264, 297) (270, 298) (275, 299) (274, 280) (228, 276) (235, 277) (234, 282) (240, 283) (284, 245) (287, 244) (288, 250) (289, 255) (292, 254) (295, 296) (300, 279) (278, 281) (285, 286) (290, 291)

\section{Algorithm for Computing Generators for Hecke Subgroups}

We can find the generators  for a representative of all the conjugacy classes of subgroups
of interest  using
 \textsf{GAP} \cite{GAP,YMR}.
  First, we use the permutation data $\sigma_{0}$, $\sigma_{1}$ obtained
from each of the Schreier coset graphs (in turn obtained from each
of the \emph{dessins}) to find the group homomorphism by images between the 
relevant Hecke group and a representative of the conjugacy class of subgroups of interest.
We then use this to define the representative in question. Finally,
we use the \textsf{GAP} command \texttt{IsomorphismFpGroup(G)}, which
returns an isomorphism from the given representative to a finitely presented
group isomorphic to that representative. This function first chooses a set of generators
of the representative and then computes a presentation in terms of these generators.

To give an example, consider the clean \emph{dessin} for the octahedron. 
We can find a set of generators as $2\times2$ matrices for
a representative of the associated class of subgroups by implementing the following code in \textsf{GAP}:

\medskip{}

\begin{singlespace}
\noindent \texttt{\small gap> f:=FreeGroup(``x'',``y'');}~\\
\texttt{\small <free group on the generators {[} x, y {]}>}~\\
\texttt{\small gap> H4:=f/{[}f.1\textasciicircum{}2,f.2\textasciicircum{}4{]};}~\\
\texttt{\small <fp group on the generators {[} x, y {]}>}\\
\texttt{\small gap>hom:=GroupHomomorphismByImages(H4,Group(}~\\
\texttt{\small (1,5)(6,11)(4,12)(2,13)(3,23)(8,14)(7,18)(10,19)(9,22)(16,24)(15,17)(20,21),}~\\
\texttt{\small (1,2,3,4)(5,6,7,8)(9,10,11,12)(13,14,15,16)(17,18,19,20)(21,22,23,24)),}~\\
\texttt{\small GeneratorsOfGroup(H4),}\\
\noindent \texttt{\small {[}(1,5)(6,11)(4,12)(2,13)(3,23)(8,14)(7,18)(10,19)(9,22)(16,24)(15,17)(20,21),}~\\
\texttt{\small (1,2,3,4)(5,6,7,8)(9,10,11,12)(13,14,15,16)(17,18,19,20)(21,22,23,24)){]});}~\\
\texttt{\small {[} x, y {]} ->}~\\
\texttt{\small {[} (1,5)(2,13)(3,23)(4,12)(6,11)(7,18)(8,14)(9,22)(10,19)(15,17)(16,24)(20,21),}~\\
\texttt{\small{} (1,2,3,4)(5,6,7,8)(9,10,11,12)(13,14,15,16)(17,18,19,20)(21,22,23,24)
{]}}~\\
\texttt{\small gap> octahedron\_group:=PreImage(hom,Stabilizer(Image(hom),1));}~\\
\texttt{\small Group(<fp, no generators known>)}~\\
\texttt{\small gap> iso:=IsomorphismFpGroup(octahedron\_group);}~\\
\texttt{\small {[} <{[} {[} 1, 1 {]} {]}|y{*}x{*}y{*}x\textasciicircum{}-1{*}y{*}x\textasciicircum{}-1>,}~\\
\texttt{\small <{[} {[} 2, 1 {]} {]}|y\textasciicircum{}-1{*}x{*}y\textasciicircum{}-1{*}x\textasciicircum{}-1{*}y\textasciicircum{}-1{*}x\textasciicircum{}-1>,}~\\
\texttt{\small <{[} {[} 3, 1 {]} {]}|y\textasciicircum{}2{*}x{*}y{*}x\textasciicircum{}-1{*}y{*}x\textasciicircum{}-1{*}y\textasciicircum{}
-1>,}~\\
\texttt{\small <{[} {[} 4, 1 {]} {]}|y\textasciicircum{}-1{*}x{*}y{*}x{*}y{*}x\textasciicircum{}-1{*}y\textasciicircum{}-2>,}~\\
\texttt{\small <{[} {[} 5, 1 {]} {]}|x{*}y\textasciicircum{}2{*}x{*}y{*}x\textasciicircum{}-1{*}y{*}x\textasciicircum{}-1{*}y\textasciicircum{}-1{*}x\textasciicircum{}-1>,}~\\
\texttt{\small <{[} {[} 6, 1 {]} {]}|x{*}y\textasciicircum{}-1{*}x{*}y{*}x{*}y{*}x\textasciicircum{}-1{*}y\textasciicircum{}-2{*}x\textasciicircum{}-1>,}~\\
\texttt{\small <{[} {[} 7, 1 {]} {]}|y{*}x{*}y\textasciicircum{}-1{*}x{*}y{*}x{*}y\textasciicircum{}-1{*}x\textasciicircum{}-1{*}y{*}x\textasciicircum{}-1{*}y\textasciicircum{}-1{*}x\textasciicircum{}-1>
{]}}~\\
\texttt{\small -> {[} F1, F2, F3, F4, F5, F6, F7 {]}}
\end{singlespace}

\medskip{}

We define the relevant Hecke group (here $H_{4}$) in the third line.
Then, the only input in each case is the permutation data $\sigma_{0}$,
$\sigma_{1}$. Once the output has been obtained, we see the generators 
(here seven:~$[F1, F2,F3, F4, F5, F6, F7]$) in the final line,
as functions of the $x$ and $y$
for the Hecke group in question.
Now, returning to the explicit matrices for these $x$ and $y$ presented in section \eqref{s:hecke}, the only thing
left to do is to multiply together the matrices and their inverses
as indicated. This will produce the generators for each representative, as required.

\end{appendix}



\bibliographystyle{plain}

\begin{thebibliography}{99}

\bibitem{FirstCachazo}
S.~Ashok, F.~Cachazo, and E.~DellÕAquila, ``Strebel Differentials with Integral Lengths and Argyres-Douglas Singularities'', arXiv:hep-th/0610080, 2006.

\bibitem{CachazoDessins}
S.~Ashok, F.~Cachazo, and E.~DellÕAquila, ``ChildrenÕs Drawings from Seiberg-Witten Curves'', arXiv:hep-th/0611082, 2006.

  \bibitem{HR2015}
  Yang-Hui He and James Read, ``Dessins d'Enfants in $\mathcal{N}=2$ Generalised Quiver Theories'', Journal of High Energy Physics 85, 2015.
  
  \bibitem{1}
 Yang-Hui He and John McKay,
  ``$\mathcal{N}=2$ Gauge Theories: Congruence Subgroups, Coset Graphs and Modular Surfaces,''
  J.\ Math.\ Phys.\  {\bf 54}, 012301, 2013.

\bibitem{He:2013lha} 
 Yang-Hui He and John McKay,
  ``Eta Products, BPS States and K3 Surfaces,''
Journal of High Energy Physics,
January 2014, 2014:113.
  
  \bibitem{YMR}
Yang-Hui He, John McKay and James Read,
``Modular Subgroups, Dessins d'Enfants and Elliptic K3 Surfaces'',
LMS Journal of Computation and Mathematics 16:271--318, 2013.
  
\bibitem{T2} 
  A.~Hanany, Y.~H.~He, V.~Jejjala, J.~Pasukonis, S.~Ramgoolam and D.~Rodriguez-Gomez,
  ``The Beta Ansatz: A Tale of Two Complex Structures'',
  JHEP {\bf 1106}, 056 (2011)
  [arXiv:1104.5490 [hep-th]].
  
  \bibitem{T1}
  V.~Jejjala, S.~Ramgoolam and D.~Rodriguez-Gomez,
  ``Toric CFTs, Permutation Triples and Belyi Pairs'',
  JHEP {\bf 1103}, 065 (2011)
  [arXiv:1012.2351 [hep-th]].
  
  \bibitem{S2}
Olivia Dumitrescu, Motohico Mulase, Brad Safnuk and Adam Sorkin, ``The spectral curve of the Eynard-Orantin recursion via the Laplace transform'', Dzhamay, Maruno and Pierce (eds.), \emph{Algebraic and Geometric Aspects of Integrable Systems and Random Matrices}, Contemporary Mathematics 593, 263-315, 2013.
  
\bibitem{S1}
Maxim Kazarian and Peter Zograf, ``Virasoro Constraints and Topological Recursion for GrothendieckÕs Dessin Counting'', Letters in Mathematical Physics, 2015.  

\bibitem{HN}
  A.~Hanany and Y.~H.~He,
  ``Non-Abelian finite gauge theories'',
  JHEP {\bf 9902}, 013 (1999)
  [hep-th/9811183].
  
  \bibitem{Zvon}
 Nicolas Magot and Alexander Zvonkin,
 ``Belyi Functions for Archimedean Solids'',
 Discrete Mathematics 217, 2000.
  
  \bibitem{DessinBook}
Ernesto Girondo and Gabino Gonzalez-Diez,
``Introduction to Compact Riemann Surfaces and Dessins d'Enfants,''
London Mathematical Society Student Texts 79, Cambridge University Press, 2012.
  
  \bibitem{4}
Alexander Zvonkin,
\newblock ``{F}unctional {C}omposition is a {G}eneralized {S}ymmetry.''
\newblock {\em Symmetry: Culture and Science}, 22:391--426, 2011.
  
  \bibitem{2}
John McKay and Abdellah Sebbar,
\newblock ``{J}-{I}nvariants of {A}rithmetic {S}emistable {E}lliptic {S}urfaces
  and {G}raphs.''
\newblock {\em Proceedings on Moonshine and related topics (Montr{\'e}al, QC,
  1999)}, pages 119--130, 2001.
  
   \bibitem{JS1}
 G.~A.~Jones and D.~Singerman, ``Theory of maps on orientable surfaces''.
 Proc. London Math. Soc. (3) 37:2, 1978, 273-307.
 
 \bibitem{JS2}
 G.~A.~Jones and D.~Singerman, ``Belyi functions, hypermaps and Galois groups''.
Bull. London Math. Soc. 28:6, 1996, 561-590.
  
  \bibitem{normal}
I.~Cangul and D. Singerman, ``Normal subgroups of Hecke groups and regular maps.''
Mathematical Proceedings of the Cambridge Philosophical Society 123, 1998.
  
  \bibitem{Jones}
 G.~A.~Jones, ``Non-congruence subgroups of the Hecke groups $G\left(\sqrt{d}\right)$''.
 Eingegangen 16, 1979.
 
  \bibitem{HeckeDessins}
Ioannis Ivrissimtzis and David Singerman, ``Regular maps and principal congruence subgroups of Hecke groups,''
European Journal of Combinatorics 26, 2005.
  
  \bibitem{HeckeSubs}
Sebahattin Ikikardes, Recep Sahin and I.~Naci Cangul,
``Principal congruence subgroups of the Hecke groups and related results'',
Bulletin of the Brazilian Mathematical Society 40, 2009.
  
  \bibitem{7}
Fritz Beukers and Hans Montanus, ``Explicit calculation of elliptic fibrations of K3-surfaces and their Belyi-maps,'' 
Number theory and polynomials, 33--51, LMS Lecture Note Series, 352, 
CUP, 2008.

 \bibitem{LMS2}
Pierre Lochak and Leila Schneps (eds.), ``Geometric Galois actions 1''. London Mathematical Society Lecture Note Series,
242, Cambridge University Press (1997).
 
 \bibitem{LMS3}
Pierre Lochak and Leila Schneps (eds.) ``Geometric Galois actions 2''. London Mathematical Society Lecture Note Series,
243, Cambridge University Press (1997).

   \bibitem{LMS1}
Leila Schneps (ed.), ``The Grothendieck theory of dessins dÕenfants''. London Mathematical Society
Lecture Note Series, 200, Cambridge University Press (1994).

  \bibitem{dessincoset}
David Singerman and J\"{u}rgen Wolfart, ``Cayley Graphs, Cori Hypermaps, and Dessins d'Enfants.''
{\em ARS Mathematica Contemporanea}. 1:144--153, 2008.
  
    \bibitem{3}
Abdellah Sebbar,
\newblock ``{M}odular {S}ubgroups, {F}orms, {C}urves and {S}urfaces.''
\newblock {\em Canad. Math. Bull.}, 45:294--308, 2002.
  
  \bibitem{Harvey}
William Harvey
``Teichm\"{u}ller Spaces, Triangle Groups and Grothendieck Dessins'',
in \emph{Handbook of Teichm\"{u}ller Theory: Volume 1}, Athanase Papadopoulos (ed.), 2007.
  
  \bibitem{G}
Branko Gr\"{u}nbaum, ``An enduring error'', Elemente der Mathematik 64 (3): 89--101, 2009.
  
  \bibitem{F}
J.~Field, ``Rediscovering the Archimedean Polyhedra:~Piero della Francesca, Luca Pacioli, Leonardo da Vinci, Albrecht D\"{u}rer, Daniele Barbaro, and Johannes Kepler'', Archive for History of Exact Sciences 50, 1997.
  
  \bibitem{Sommerville}
D.~M.~Y.~Sommerville, ``Semi-Regular Networks of the Plane in Absolute Geometry'', Transactions of the Royal Society of Edinburgh 41, 1905.
  
  \bibitem{BPR}
Egbert Brieskorn, Anna Pratoussevitch, Frank Rothenh\"ausler,
``The combinatorial geometry of singularities and Arnold's series E, Z, Q'',
Moscow Mathematical Journal, v.3, 273-333, 741 (2003).
  
  \bibitem{GAP}
  The GAP Group, \emph{GAP -- Groups, Algorithms, and Programming, 
  Version 4.5.6}; 
  2012, \verb+(http://www.gap-system.org)+.
  
  \bibitem{sage}
  The Sage Group, \emph{Sage -- Version 5.3};
  2012,
  \verb+(http://www.sagemath.org/)+.

\bibitem{Gaiotto}
Davide Gaiotto,
``$\mathcal{N} = 2$ Dualities,'' Journal of High Energy Physics, vol. 34, no. 8, 2012.

\bibitem{Dummies}
Yuji Tachikawa, ``$\mathcal{N}=2$ Supersymmetric Dynamics for Pedestrians'', 2014, arXiv:1312.2684v2.



\bibitem{ClarkVoight}
Pete Clark and John Voight,
``Algebraic Curves Uniformized by Congruence Subgroups of Triangle Groups'', 2011.

 \bibitem{PrinCon}
 M.~L.~Lang, C.~H.~Lim, S. P. Tan, ``Principal congruence subgroups of the Hecke
groups''. J. Number theory, vol. 85, issue 2 (2000), 220-230.












\end{thebibliography}

\end{document}